\documentclass[12pt, twoside, eqno]{article}
\usepackage{latexsym}
\usepackage{amssymb}
\usepackage{amsfonts}
\usepackage{amsmath}
\begin{document}

\noindent {\bf \large Fuzzy right (left) ideals in hypergroupoids and 
fuzzy bi-ideals in hypersemigroups}
\bigskip

\noindent{\bf Niovi Kehayopulu}\\
May 25, 2016\bigskip

\smallskip{\small

\noindent{\bf Abstract.} We introduce the concepts of fuzzy right and 
fuzzy left ideals of hypergroupoids and the concept of a fuzzy 
bi-ideal of an hypersemigroup and we show that a fuzzy subset $f$ of 
an hypergroupoid $H$ is a fuzzy right (resp. fuzzy left) ideal of $H$ 
if and only if $f\circ 1\preceq f$ (resp. $1\circ f\preceq f)$ and 
for an hypersemigroup $H$, a fuzzy subset $f$ of $H$ is a bi-ideal of 
$H$ if and only if $f\circ 1\circ f\preceq f$. These 
characterizations are very useful for the investigation. The paper 
serves as an example to show the way we pass from fuzzy groupoids 
(semigroups) to fuzzy hypergroupoids (hypersemigroups).\medskip

}
\section{Introduction} In our paper in [2] we gave, among others, 
some equivalent definitions of fuzzy right (left) ideals and fuzzy 
bi-ideals in ordered semigroups which are very useful for 
applications. Besides, they show how similar is the theory of ordered 
semigroups based on fuzzy ideals with the theory of ordered 
semigroups based on ideals or on ideal elements. Using these 
definitions many results on fuzzy ordered semigroups or on fuzzy 
semigroups (without order) can be drastically simplified. The present 
paper is based on our paper in [2], and the aim is to show the way we 
pass from fuzzy groupoids (semigroups) to fuzzy hypergroupoids 
(hypersemigroups).
\section{Main results}\noindent An {\it hypergroupoid} is a nonempty 
set $H$ with an hyperoperation $$\circ : H\times H \rightarrow {\cal 
P}^*(H) \mid (a,b) \rightarrow a\circ b$$on $H$ and an operation $$* 
: {\cal P}^*(H)\times {\cal P}^*(H) \rightarrow {\cal P}^*(H) \mid 
(A,B) \rightarrow A*B$$ on ${\cal P}^*(H)$ (induced by the operation 
of $H$) such that $$A*B=\bigcup\limits_{(a,b) \in\,A\times B} 
{(a\circ b)}$$ for every $A,B\in {\cal P}^*(H)$ (${\cal P}^*(H)$ 
denotes the set of nonempty subsets of $H$).\smallskip

The operation ``$*$" is well defined. Indeed: If $(A,B)\in {\cal 
P}^*(H) \times {\cal P}^*(H)$, then $A*B=\bigcup\limits_{(a,b) 
\in\,A\times B} {(a\circ b)}$. For every $(a,b)\in A\times B$, we 
have $(a,b)\in H\times H$, then $(a\circ b)\in {\cal P}^*(H)$, thus 
we get $A*B\in {\cal P}^*(H)$. If $(A,B),(C,D)\in {\cal P}^*(H)\times 
{\cal P}^*(H)$ such that $(A,B)=(C,D)$, then 
$$A*B=\bigcup\limits_{(a,b) \in\,A\times B} {(a\circ 
b)}=\bigcup\limits_{(a,b) \in\,C\times D} {(a\circ b)}=C*D.$$
As the operation ``$*$" depends on the hyperoperation ``$\circ$", an 
hypergroupoid can be also denoted by $(H,\circ)$ (instead of 
$(H,\circ,*)$).

If $H$ is an hypergroupoid then, for any $x,y\in H$, we have $x\circ 
y=\{x\}*\{y\}$. Indeed,$$\{x\}*\{y\}=\bigcup\limits_{u \in \{ x\} ,v 
\in \{ y\} } {u \circ v = x \circ y}.$$

An hypergroupoid $H$ is called {\it hypersemigroup} if $$(x\circ 
y)*\{z\}=\{x\}*(y\circ z)$$ for every $x,y,z\in H$. Since $x\circ 
y=\{x\}*\{y\}$ for any $x,y\in H$, an hypergroupoid $H$ is an 
hypersemigroup if and only if, for any $x,y,z\in H$, we 
have$${\Big(}\{x\}*\{y\}{\Big)}*\{z\}=\{x\}*{\Big(}\{y\}*\{z\}{\Big)}.$$

Following Zadeh, if $(H,\circ)$ is an hypergroupoid, we say that $f$ 
is a fuzzy subset of $H$ (or a fuzzy set in $H$) if $f$ is a mapping 
of $H$ into the real closed interval $[0,1]$ of real numbers, that is 
$f : H \rightarrow [0,1]$. For an element $a$ of $H$, we denote by 
$A_a$ the subset of $H\times H$ defined as follows:$$A_a:=\{(y,z)\in 
H\times H \mid a\in y\circ z\}.$$For two fuzzy subsets $f$ and $g$ of 
$H$, we denote by $f\circ g$ the fuzzy subset of $H$ defined as 
follows:
$$f \circ g: H \to [0,1]\,\,\,a \to \left\{ \begin{array}{l}
\bigvee\limits_{(y,z) \in {A_a}} {\min \{ f(y),g(z)\} 
\,\,\,\,if\,\,\,{A_a} \ne \emptyset } \\
\,\,\,\,0\,\,\,\,\,\,if\,\,\,\,{A_a} = \emptyset. \,\,\,\,\,\,\,\,\,
\end{array} \right.$$Denote by $F(H)$ the set of all fuzzy subsets of 
$H$ and by ``$\preceq$" the order relation on $F(H)$ defined 
by:$$f\preceq g \;\Longleftrightarrow\; f(x)\le g(x) \mbox { for 
every } x\in H.$$We finally show by 1 the fuzzy subset of $H$ defined 
by:$$1: H \rightarrow [0,1] \mid x \rightarrow 1(x):=1.$$Clearly, the 
fuzzy subset 1 is the greatest element of the ordered set 
$(F(H),\preceq)$ (that is, $1\succeq f$ $\forall f\in F(H)$).

For two fuzzy subsets $f$ and $g$ of an hypergroupoid $H$ we denote 
by $f\wedge g$ the fuzzy of $H$ defined as follows:$$f\wedge g : 
H\rightarrow [0,1] \mid x\rightarrow (f\wedge 
g)(x):=\min\{f(x),g(x)\}.$$One can easily prove that the fuzzy subset 
$f\wedge g$ is the infimum of the fuzzy subsets $f$ and $g$, so we 
write $f\wedge g=\inf\{f,g\}$.

We denote the hyperoperation on $H$ and the multiplication between 
the two fuzzy subsets of $H$ by the same symbol (no confusion is 
possible).\medskip

The following proposition, though clear, plays an essential role in 
the theory of hypergroupoids.\medskip

\noindent{\bf Proposition 1.} {\it Let $(H,\circ)$ be an 
hypergroupoid, $x\in H$ and $A,B\in {\cal P}^*(H)$. Then we have the 
following:

$(1)$ $x\in A*B$ $\Longleftrightarrow$ $x\in a\circ b$ for some $a\in 
A$, $b\in B$.

$(2)$ If $a\in A$ and $b\in B$, then $a\circ b\subseteq A*B$.} 
\medskip

\noindent{\bf Lemma 2.} {\it Let H be an hypergroupoid. Then we have 
the following:

$(1)$ If $A\subseteq B$ and $C\subseteq D$, then $A*C\subseteq B*D$ 
and $C*A\subseteq D*B$.

$(2)$ $H*H\subseteq H$.}\medskip

\noindent{\bf Proof.} (1) Let $A\subseteq B$, $C\subseteq D$ and 
$x\in A*C$. By Proposition 1(1), there exist $a\in A$ and $c\in C$ 
such that $x\in a\circ c$. Since $a\in B$ and $c\in D$, by 
Proposition 1(2), we have $a\circ c\in B*D$, then $x\in B*D$. 
Similarly, $C*A\subseteq D*B$.\smallskip

\noindent(2) $H*H\subseteq H$. Indeed: $H*H=\bigcup\limits_{u\in H, 
\, v\in H} {u\circ v}$. On the other hand, for any $u,v\in H$, we 
have $u\circ v\subseteq H$. Thus we have $H*H\subseteq H$. 
$\hfill\Box$\medskip

\noindent{\bf Definition 3.} Let $H$ be an hypergroupoid. A fuzzy 
subset $f$ of $H$ is called a {\it fuzzy right ideal} of $H$ 
if$$f(x\circ y)\ge f(x) \mbox { for every } x,y\in H,$$in the sense 
that: if $x,y\in H$ and $u\in x\circ y$, then $f(u)\ge 
f(x)$.\medskip

\noindent{\bf Theorem 4.} {\it Let H be an hypergroupoid  and f a 
fuzzy subset of H. Then f is a fuzzy right ideal of H if and only if 
$$f\circ 1\preceq f.$$}{\bf Proof.} $\Longrightarrow$. Let $a\in H$. 
Then $(f\circ 1)(a)\le f(a)$. In fact: \\
(i) If $A_a=\emptyset$, then $(f\circ 1)(a):=0$. Since $f$ is a fuzzy 
subset of $H$, we have $f(a)\ge 0$. Thus we have $(f\circ 1)(a)\le 
f(a)$.\\
(ii) Let $A_a\not=\emptyset$. Then\begin{equation}(f \circ 1)(a): = 
\bigvee\limits_{(y,z) \in {A_a}} {\min \{ f(y),1(z)\} }= 
\bigvee\limits_{(y,z) \in {A_a}}f(y) \tag{$*$}\end{equation}On the 
other hand, $f(y)\le f(a) \mbox { for every } (y,z)\in A_a 
$\hfill$(**)$\\Indeed, if $(y,z)\in A_a$, then $a\in y\circ z$ and, 
since $f$ is a fuzzy right ideal of $H$, we have $f(a)\ge f(y)$.

By $(**)$, we have $\bigvee\limits_{(y,z) \in {A_a}}f(y) \le f(a).$ 
Then, by $(*)$, $(f\circ 1)(a)\le f(a)$.\smallskip

\noindent$\Longleftarrow$. Let $x,y\in H$ and $u\in x\circ y$. Then 
$f(u)\ge f(x)$. Indeed: \\Since $x,y\in H$ and $u\in x\circ y$, we 
have $(x,y)\in A_u$.
Since $A_u\not=\emptyset$, we have\begin{eqnarray*}(f\circ 
1)(u):=\bigvee\limits_{(t,s) \in {A_u}} {\min \{ f(t),1(s)\} 
}=\bigvee\limits_{(t,s) \in {A_u}}f(t)\ge f(t)
\;\; \forall \;(t,s)\in A_u. \end{eqnarray*}Since $(x,y)\in A_u$, we 
have $(f\circ 1)(u)\ge f(x)$. Since $u\in x\circ y\subseteq H$, we 
have $u\in H$. Since $u\in H$, by hypothesis, we have $(f\circ 
1)(u)\le f(u)$. Then we have $f(u)\ge f(x)$. $\hfill\Box$\medskip

\noindent{\bf Definition 5.} Let $H$ be an hypergroupoid. A fuzzy 
subset $f$ of $H$ is called a {\it fuzzy left ideal} of $H$ 
if$$f(x\circ y)\ge f(y) \mbox { for every } x,y\in H,$$in the sense 
that: if $x,y\in H$ and $u\in x\circ y$, then $f(u)\ge f(y)$.\\In a 
similar way as in Theorem 4, we can prove the following 
theorem.\medskip

\noindent{\bf Theorem 6.} {\it A fuzzy subset f of an hypergroupoid H 
is a fuzzy left ideal of H if and only if, for any fuzzy subset $f$ 
of $H$, we have$$1\circ f\preceq f.$$}
\noindent{\bf Definition 7.} Let $H$ be an hypergroupoid. A fuzzy 
subset $f$ of $H$ is called a {\it fuzzy quasi-ideal} of $H$ if$$x\in 
b\circ s \mbox { and } x\in t\circ c\; \Longrightarrow\; f(x)\ge 
\min\{f(b), f(c)\} \;\forall\;x,b,s,t,c\in H.$${\bf Theorem 8.} {\it 
Let H be an hypergroupoid. A fuzzy subset f of H is a fuzzy 
quasi-ideal of H if and only if$$(f\circ 1)\wedge (1\circ f)\preceq 
f.$$}{\bf Proof.} $\Longrightarrow$. Let $x\in H$. Then 
${\Big(}(f\circ 1)\wedge (1\circ f){\Big)}(x)\le f(x)$, that is\\
$\min\{(f\circ 1)(x),(1\circ f)(x)\}\le f(x)$. Indeed: For 
$A_x=\emptyset$ this is clear.\\Let $A_x\not=\emptyset$. Then
$$(f \circ 1)(x): = \mathop  \bigvee \limits_{(y,s) \in {A_x}} \min 
\{ f(y),1(s)\}  = \mathop  \bigvee \limits_{(y,s) \in {A_x}} f(y) 
$$and$$(1 \circ f)(x): = \mathop  \bigvee \limits_{(t,z) \in {A_x}} 
\min \{ 1(t),f(z)\}  = \mathop  \bigvee \limits_{(t,z) \in {A_x}} 
f(z).$$If $f(x)\ge (f\circ 1)(x)$, then $f(x)\ge \min\{(f\circ 1)(x), 
(1\circ f)(x)\}$.\\Let $f(x)<(f\circ 1)(x)$. Then there exists 
$(y,s)\in A_x$ such that $f(y)>f(x) \;\;\;\;\;\;\;\;\;\; (*)$ 
(otherwise $(f\circ 1)(x)\le f(x)$ which is impossible).\\We prove 
that $f(x)\ge f(z)$ for every $(t,z)\in A_x$. Then we have$$f(x)\ge 
\mathop  \bigvee \limits_{(t,z) \in {A_x}} f(z) = (1 \circ f)(x)\ge 
\min \{ (f \circ 1)(x),(1 \circ f)(x)\}$$and the proof is complete.

Let now $(t,z)\in A_x$. Then $f(x)\ge f(z)$. Indeed: Since $(y,s)\in 
A_x$, we have $y,s\in H$ and $x\in y\circ s$. Since $(t,z)\in A_x$, 
we have $t,z\in H$ and $x\in t\circ z$. Since $x,y,s,t,z\in H$ such 
that $x\in y\circ s$ and $x\in t\circ z$, by hypothesis, we have 
$f(x)\ge\min\{f(y),f(z)\}$. If $\min\{f(y),f(z)\}=f(y)$, then 
$f(x)\ge f(y)$ which is impossible by $(*)$. Thus we have 
$\min\{f(y),f(z)\}=f(z)$, and $f(x)\ge f(z)$. \smallskip

\noindent$\Longleftarrow$. Let $x,b,s,t,c\in H$ such that $x\in 
b\circ s$ and $x\in t\circ c$. Then $f(x)\ge\min\{f(b),f(c)\}$. 
Indeed: By hypothesis, we have$$f(x)\ge {\Big(}(f\circ 1)\wedge 
(1\circ f){\Big)}(x):=\min\{(f\circ 1)(x),(1\circ f)(x)\}.$$Since 
$x\in b\circ s$, we have $(b,s)\in A_x$, then
\[\begin{array}{l}
(f \circ 1)(x): = \bigvee\limits_{(u,v) \in {A_x}} {\min \{ 
f(u),1(v)\} }  = \bigvee\limits_{(u,v) \in {A_x}} {f(u)\ge f(b).} 
\end{array}\]Since $x\in t\circ c$, we have $(t,c)\in A_x$, then
\[\begin{array}{l}
(1 \circ f)(x): = \bigvee\limits_{(w,k) \in {A_x}} {\min \{ 
1(w),f(k)\} }  = \bigvee\limits_{(w,k) \in {A_x}} {f(k)\ge f(c).} 
\end{array}\]Thus we have$$f(x)\ge \min\{(f\circ 1)(x),(1\circ 
f)(x)\}\ge \min\{f(b),f(c)\}.$$ $\hfill\Box$

The proof of the associativity of fuzzy sets on semigroups given in 
[1] can be naturally transferred to hypersemigroups in the following 
proposition.\medskip

\noindent{\bf Proposition 9.} {\it If $H$ is an hypersemigroup, then 
the set of all fuzzy subsets of $H$ is a semigroup.}\medskip

\noindent{\bf Proof.} Let $f, g, h$ be fuzzy subsets of $H$ and $a\in 
H$. Then$${\Big(}(f\circ g)\circ h{\Big)}(a)={\Big(}f\circ (g\circ 
h){\Big)}(a).$$Indeed: If $A_a=\emptyset$, then ${\Big(}(f\circ 
g)\circ h{\Big)}(a)=0={\Big(}f\circ (g\circ h){\Big)}(a).$

Let $A_a\not=\emptyset$. Then$${\Big(}(f\circ g)\circ h{\Big)}(a) = 
\bigvee\limits_{(y,z)\in A_a} {\min \{ (f \circ g)(y),h(z)\} } 
$$and$${\Big(}f\circ (g\circ h{\Big)}(a) = \bigvee\limits_{(u,v)\in 
A_a} {\min \{ f(u),(g\circ h)(v)\}}.$$We put 
\begin{equation}t:=\bigvee\limits_{(y,z)\in A_a} {\min \{ (f \circ 
g)(y),h(z)\} } \tag{$*$}\end{equation}$$s:=\bigvee\limits_{(u,v)\in 
A_a} {\min \{ f(u),(g\circ h)(v)\}}.$$We prove that $t\ge 
\min\{(f(u),(g\circ h)(v)\}$ for every $(u,v)\in A_a$. Then we have 
$t\ge s$. In a similar way we prove that $s\ge t$, and so 
$s=t$.\smallskip

Let $(u,v)\in A_a$ ($\Rightarrow$ $t\ge \min\{f(u),(g\circ 
h)(v)\}\;\; ?\,)$\\
(A) Let $t\ge f(u)$. Since $f(u)\ge \min\{f(u),(g\circ h)(v)\}$, we 
have $$t\ge \min\{f(u),(g\circ h)(v)\}.$$(B) Let $t<f(u)$. We 
consider the cases:

(a) Let $A_v=\emptyset$. Then $(g\circ h)(v):=0$. Since $f$ is a 
fuzzy set in $H$, we have $f(u)\ge 0$, then $\min\{f(u),(g\circ 
h)(v)\}=0$. Since $t\in [0,1]$, we have $t\ge 0$, so $t\ge 
\min\{f(u),(g\circ h)(v)\}$.

(b) Let $A_v\not=\emptyset$. Then $(g \circ h)(v)= 
\bigvee\limits_{(c,d)\in A_v} {\min \{ g({c}} ),h({d})\}$. We prove 
that$$t\ge\min\{g(c),h(d)\}\mbox { for every } (c,d)\in A_v.$$ Then 
we have $t \ge (g \circ h)(v)\ge \min\{f(u),(g\circ 
h)(v)\}.$\smallskip

Let now $(c,d)\in A_v$ ($\Rightarrow$ $t\ge 
\min\{g(c),h(d)\}\;\;?\,)$

\noindent (i) Let $t\ge g(c)$. Since $g(c)\ge\min\{g(c),h(d)\}$, we 
have $t\ge \min\{g(c),h(d)\}$.

\noindent (ii) Let $t<g(c)$. Since $(u,v)\in A_a$, we have $a\in 
u\circ v$. Since $(c,d)\in A_v$, we have and $v\in c\circ d$. Then
we have\begin{eqnarray*}a\in u\circ v&=&\{u\}*\{v\}\subseteq 
\{u\}*(c\circ d) \mbox { (by Lemma 2)}\\&=&(u\circ c)*\{d\} \mbox { 
(the operation } ``*" \mbox { is associative)}.\end{eqnarray*}By 
Proposition 1(1), $a\in w\circ d$ for some $w\in u\circ c$. Since 
$a\in w\circ d$, we have $(w,d)\in A_a$ then, by $(*)$, we 
have\begin{equation}t\ge \min\{(f\circ g)(w),h(d)\} 
\tag{$**$}\end{equation}Since $w\in u\circ c$, we have $(u,c)\in 
A_w$, then\begin{equation}(f\circ g)(w)= \bigvee\limits_{(l,k)\in 
A_w} {\min \{ f({l}} ),g({k})\}\ge \min\{f(u),g(c)\} 
\tag{$***$}\end{equation}Since $t<f(u)$ and $t<g(c)$, we have 
$t<\min\{f(u),g(c)\}$. Then, by (***), $t<(f\circ g)(w)$. Then, by 
$(**)$, $t\ge h(d)$. Since $h(d)\ge \min\{g(c),h(d)\}$, we have $t\ge 
\min\{g(c),h(d)\}$ and the proof of the proposition is complete. 
$\hfill\Box$\medskip

According to Proposition 9, for any fuzzy subsets $f,g,h$ of $H$, we 
write

$(f\circ g)\circ h=f\circ (g\circ h):=f\circ g\circ h$.\medskip

\noindent{\bf Definition 10.} Let $H$ be an hypersemigroup. A fuzzy 
subset $f$ of $H$ is called a {\it fuzzy bi-ideal} of $H$ 
if$$f{\Big(}(x\circ y)*\{z\}{\Big)}\ge \min \{f(x),f(z)\} 
\;\;\forall\; x,y,z\in H,$$in the sense that if $u\in (x\circ 
y)*\{z\}$, then $f(u)\ge \min \{f(x),f(z)\}$.\medskip

\noindent{\bf Theorem 11.}  {\it Let H be an hypersemigroup and f a 
fuzzy subset of S. Then f is a fuzzy bi-ideal of H if and only 
if$$f\circ 1\circ f\preceq f.$$}{\bf Proof.} $\Longrightarrow$. Let 
$a\in H$. Then $(f\circ 1\circ f)(a)\le f(a)$. In fact: If 
$A_a=\emptyset$, then$$(f\circ 1\circ f)(a)={\Big(}(f\circ 1)\circ 
f{\Big)}(a):=0\le f(a).$$Let $A_a\not=\emptyset$. Then\[(f \circ 1 
\circ f)(a): = \bigvee\limits_{(y,z) \in {A_a}} {\min \{ (f \circ 
1)(y),f(z)\}. } \] It is enough to prove 
that\begin{equation}\min\{(f\circ 1)(y),f(z)\}\le f(a) \mbox { for 
every } (y,z)\in A_a \tag{$*$}\end{equation}For this purpose, let 
$(y,z)\in A_a$. If $A_y=\emptyset$, then $(f\circ 1)(y):=0\le f(z)$, 
and\\$\min\{(f\circ 1)(y),f(z)\}=0\le f(a)$. Let now 
$A_y\not=\emptyset$. Then\begin{equation}(f \circ 1)(y): = 
\bigvee\limits_{(x,w) \in {A_y}} {\min \{ 
f(x),1(w)\}=\bigvee\limits_{(x,w) \in 
{A_y}}f(x)}\tag{1}\end{equation}We consider the following cases:\\
(i) Let $f(a)\ge (f\circ 1)(y)$. Then$$f(a)\ge (f\circ 1)(y)\ge 
\min\{(f\circ 1)(y),f(z)\},$$and condition $(*)$ is satisfied.\\
(ii) Let $f(a)<(f\circ 1)(y)$. Then there exists $(x,w)\in A_y$ such 
that $f(a)< f(x) \;\, (2)$. Indeed: If $f(a)\ge f(x)$ for every 
$(x,w)\in A_y$, then $f(a)\ge \bigvee\limits_{(x,w) \in {A_y}}f(x)$. 
Then, by (1), $f(a)\ge (f\circ 1)(y)$ which is impossible.

Since $(y,z)\in A_a$, we have $y,z\in H$ and $a\in y\circ z$. Since 
$(x,w)\in A_y$, we have $x,w\in H$ and $y\in x\circ w$. Then we 
have$$a\in y\circ z=\{y\}*\{z\}
\subseteq (x\circ w)*\{z\}$$and, since $f$ is a fuzzy bi-ideal of 
$H$, we have\begin{equation}f(a)\ge \min\{f(x),f(z)\} 
\tag{3}\end{equation}If $f(x)\le f(z)$ then $f(a)\ge f(x)$, which is 
impossible by (2). Thus we 
have\begin{equation}f(x)>f(z)\tag{4}\end{equation}Then 
$\min\{f(x),f(z)\}=f(z)$ and, by (3), $f(a)\ge 
f(z)\hfill{(5)}$\\Since $(x,w)\in A_y$, by (1), we have $f(x)\le 
(f\circ 1)(y)\hfill{(6)}$\\By (4) and (6), we have $f(z)<f(x)\le 
(f\circ 1)(y)$, then$$\min\{f\circ 1)(y),f(z)\}=f(z)\le f(a) \mbox { 
(by (5)}),$$ and condition $(*)$ is satisfied.\smallskip

\noindent$\Longleftarrow$. Let $x,y,z\in H$ and $u\in (x\circ 
y)*\{z\}$. Then $f(u)\ge \min\{f(x),f(z)\}$.\\Indeed: Since $u\in 
(x\circ y)*\{z\}$, by Proposition 1(1), there exist $a\in x\circ y$ 
such that $u\in a\circ z$. Since $a\in x\circ y$, we have $(x,y)\in 
A_a$. Since $u\in a\circ z$, we have $(a,z)\in A_u$. Since $(a,z)\in 
A_u$, $A_u$ is a nonempty set and we also have$$(f \circ 1 \circ 
f)(u) = \bigvee\limits_{(t,s) \in {A_u}} {\min \{ (f \circ 
1)(t),f(s)\}  \ge \min \{ (f \circ 1)(a),f(z)\} .}$$Since $(x,y)\in 
A_a$, we have $A_a\not=\emptyset$ and we also have$$(f \circ 1)(a) = 
\bigvee\limits_{(z,w) \in {A_a}} {\min \{ f(z),1(w)\}  \ge \min \{ 
f(x),1(y)\}  = f(x).}$$Then we have $$(f\circ 1\circ f)(u)\ge 
\min\{f(x),f(z)\}.$$On the other hand, since $u\in (x\circ 
y)*\{z\}\subseteq H*H\subseteq H$ and $f\circ 1\circ f\preceq f$, we 
have $(f\circ 1\circ f)(u)\le f(u)$. Thus we have$$f(u)\ge (f\circ 
1\circ f)(u)\ge \min\{f(x),f(z)\},$$ and the proof of the theorem is 
complete. $\hfill\Box$\bigskip

\noindent{\bf Note.} The characterization of fuzzy right (left) and 
fuzzy bi-ideals of an hypersemigroup $H$ using the greatest element 
$1$ of the ordered set of fuzzy subsets of $H$ is very useful for 
further investigation. Using these definitions many proofs can be 
drastically simplified. As an example, one can immediately see that 
every fuzzy right (or fuzzy left) ideal of $H$ is a fuzzy bi-ideal of 
$H$. Indeed, if $f$ is a fuzzy right ideal of $H$, then $f\circ 
1\circ f=f\circ (1\circ f)\preceq f\circ 1\preceq f$.
{\small\bigskip

\smallskip

\noindent University of Athens, Department of Mathematics, 15784 
Panepistimiopolis, Greece\\email: nkehayop@math.uoa.gr
\end{document}